\documentclass[12pt]{article}
\usepackage{amssymb}
\usepackage{latexsym}

\newtheorem{definition}{Definition}[section]
\newtheorem{theorem}[definition]{Theorem}
\newtheorem{lemma}[definition]{Lemma}
\newtheorem{corollary}[definition]{Corollary}
\newtheorem{remark}[definition]{Remark}
\newtheorem{example}[definition]{Example}
\newtheorem{conjecture}[definition]{Conjecture}
\newtheorem{problem}[definition]{Problem}
\newtheorem{note}[definition]{Note}
\newtheorem{assumption}[definition]{Assumption}
\newtheorem{proposition}[definition]{Proposition}

\def\I{\mathbb I}

\def\K{\mathbb K}

\def\Z{\mathbb Z}
\def\K{\mathbb K}

\begin{document}

\title{\bf The $q$-tetrahedron algebra and its finite dimensional
irreducible modules}
\author{
Tatsuro Ito{\footnote{
Department of Computational Science,
Faculty of Science,
Kanazawa University,
Kakuma-machi,
Kanazawa 920-1192, Japan
}}
$\;$ and
Paul Terwilliger{\footnote{
Department of Mathematics, University of
Wisconsin, 480 Lincoln Drive, Madison WI 53706-1388 USA}
}}
\date{}
%to get date printout, comment out above line

\maketitle
\begin{abstract}
Recently B. Hartwig and the second author found
a presentation for the three-point $\mathfrak{sl}_2$ loop algebra
via generators and relations.
To obtain this presentation they defined an algebra $\boxtimes$
by generators and relations, and displayed an isomorphism from
$\boxtimes$ to the three-point $\mathfrak{sl}_2$ loop algebra.
We introduce a quantum analog of
$\boxtimes$  which we call $\boxtimes_q$.
We define $\boxtimes_q$ via generators and relations. 
We show how 
$\boxtimes_q$ is related to the quantum group
$U_q(\mathfrak{sl}_2)$,
the 
$U_q(\mathfrak{sl}_2)$ loop algebra,
and the positive part of
$U_q({\widehat{\mathfrak{sl}}_2})$.
We describe the finite dimensional irreducible
$\boxtimes_q$-modules under the assumption that
$q$ is not a root of 1, and the underlying
field is algebraically closed.
\medskip

 \end{abstract}

\section{Introduction}
In \cite{HT}
B. Hartwig and the second author
gave a presentation of the three-point $\mathfrak{sl}_2$ loop
algebra via generators and relations. To obtain this
presentation they defined a Lie algebra $\boxtimes$ by
generators and relations, and displayed an isomorphism
 from $\boxtimes$ to
the three-point $\mathfrak{sl}_2$ loop algebra.
$\boxtimes$ has essentially six generators, and it is natural
to identify these with the six edges of a tetrahedron.
For each face of the
tetrahedron the three 
surrounding edges
form a basis for a subalgebra of $\boxtimes$
that is isomorphic to $\mathfrak{sl}_2$ \cite[Corollary 12.4]{HT}. 
Any five of the six edges of the tetrahedron generate a subalgebra of
$\boxtimes$ 
that is isomorphic to the $\mathfrak{sl}_2$ loop algebra
 \cite[Corollary 12.6]{HT}.
Each pair of opposite
edges generate a subalgebra of $\boxtimes$ that is isomorphic
to the Onsager algebra \cite[Corollary 12.5]{HT}.
Let us call these Onsager subalgebras. Then $\boxtimes
$
is the direct sum of its three Onsager subalgebras
 \cite[Theorem 11.6]{HT}.
%The universal central extension of 
%of three-point $\mathfrak{sl}_2$ loop algebra
%was obtained Bremner \cite{Br}.
%A presentation for this
%extension by generators and relations was obtained by Benkart and the second
%author \cite{BT},  
%The irreducible finite dimensional $\boxtimes$-modules
%were obtained by Hartwig \cite{Ha}.
In \cite{E} Elduque found an 
attractive decomposition of $\boxtimes$ into a direct sum of three 
abelian subalgebras, and he
showed how these subalgebras 
are related to the Onsager subalgebras. 
In \cite{PT} Pascasio and the second author give an
action of $\boxtimes$ on the standard module of the Hamming association
scheme. In \cite{Br} Bremner
obtained 
the universal central extension of 
the three-point $\mathfrak{sl}_2$ loop algebra.
By modifying the defining relations for $\boxtimes$, 
Benkart and the second author obtained  
a presentation for this
extension by generators and relations
 \cite{BT}. 
In \cite{Ha} Hartwig obtained the 
irreducible finite dimensional $\boxtimes$-modules
over an algebraically closed field with characteristic 0.

\medskip
\noindent 
In this paper we introduce a quantum analog of
$\boxtimes$ which we call 
$\boxtimes_q$.
We define $\boxtimes_q$ using generators and relations.
We show how $\boxtimes_q$ is related to the
quantum group
$U_q(\mathfrak{sl}_2)$ in roughly the same
way that $\boxtimes $ is related to $\mathfrak{sl}_2$.
We show how $\boxtimes_q$ is related to
the 
$U_q(\mathfrak{sl}_2)$ loop algebra
in roughly the same way that
$\boxtimes$ is related to the $\mathfrak{sl}_2$ loop algebra.
In \cite{NN} we considered an algebra
${\mathcal A}_q$ on two generators
subject to the cubic $q$-Serre relations.
${\mathcal A}_q$ is often called the
{\it positive part of 
$U_q({\widehat{\mathfrak{sl}}_2})$}.
We show how $\boxtimes_q$ is related to
${\mathcal A}_q$ in roughly the same way that
$\boxtimes$ is related to the Onsager algebra.
We describe the finite dimensional irreducible
$\boxtimes_q$-modules under the assumption
 that $q$ is not a root of 1, and the
 underlying field is algebraically closed.
As part of our description we relate
these $\boxtimes_q$-modules to the following
kind of ${\mathcal A}_q$-module.
 A finite dimensional irreducible 
${\mathcal A}_q$-module $V$ is called {\it NonNil} whenever
neither of the two  ${\mathcal A}_q$ generators is nilpotent
on $V$. We described these ${\mathcal A}_q$-modules in \cite{NN}.
Associated with such a module there 
is an ordered pair of nonzero scalars
called its {\it type}. Also, for each finite dimensional irreducible
$\boxtimes_q$-module there is a scalar in the set $\lbrace 1,-1\rbrace$
called its {\it type}.
The main result of this paper is an explicit bijection between the following
two sets: 
\begin{enumerate}
\item[{\rm (i)}] 
 the isomorphism classes of finite dimensional
 irreducible $\boxtimes_q$-modules of type $1$;
\item[{\rm (ii)}]  
the isomorphism classes of 
NonNil finite dimensional irreducible 
 ${\mathcal A}_q$-modules of type $(1,1)$.
\end{enumerate}
\noindent
All of the original results in this paper
are about $\boxtimes_q$, although we
will initially discuss $\boxtimes$
in order to motivate things.
The paper is organized as follows.
In Section 2 we define  $\boxtimes$
and discuss how it is related to $\mathfrak{sl}_2$.
In Section 3 we recall how $\boxtimes$ is related
to the $\mathfrak{sl}_2$ loop algebra. In Section
4 we recall how $\boxtimes$ is related to the
Onsager algebra. In Section 5 we discuss the
finite dimensional irreducible $\boxtimes$-modules.
In Section 6 we introduce the algebra $\boxtimes_q$.
In Section 7 we discuss how $\boxtimes_q$ is related
to $U_q(\mathfrak{sl}_2)$.
In Section 8 we discuss how $\boxtimes_q$ is related
to the $U_q(\mathfrak{sl}_2)$ loop algebra.
In Section 9 we discuss how $\boxtimes_q$ is related
to the algebra ${\mathcal A}_q$.
In Sections 10--18 we describe the finite dimensional
irreducible $\boxtimes_q$-modules.
In Section 19 we mention some open problems.

\medskip
\noindent Throughout the paper $\K$ denotes a field.
\section{The tetrahedron algebra}

In this section we recall the tetrahedron algebra
and discuss how it is related to
$\mathfrak{sl}_2$.

\medskip
\noindent
Until further notice  assume
the field $\K$ has characteristic 0.

\begin{definition}
\label{def:tet}
\rm
\cite{HT}
Let $\boxtimes$ denote the Lie algebra over $\K$
that has  generators
\begin{eqnarray}
\label{eq:boxgen}
\lbrace x_{ij} \,|\,i,j\in \I, i\not=j\rbrace
\qquad \qquad \I = \lbrace 0,1,2,3\rbrace
\end{eqnarray}
and the following relations:
\begin{enumerate}
\item[{\rm (i)}]  
For distinct $i,j\in \I$,
\begin{eqnarray*}
x_{ij}+x_{ji} = 0.
\label{eq:rel0}
\end{eqnarray*}
\item[{\rm (ii)}]  
For mutually distinct $h,i,j\in \I$,
\begin{eqnarray*}
\lbrack x_{hi},x_{ij}\rbrack = 2x_{hi}+2x_{ij}.
\label{eq:rel1}
\end{eqnarray*}
\item[{\rm (iii)}]  
For mutually distinct $h,i,j,k\in \I$,
\begin{eqnarray*}
\lbrack x_{hi},
\lbrack x_{hi},
\lbrack x_{hi},
x_{jk}\rbrack \rbrack \rbrack= 
4 \lbrack x_{hi},
x_{jk}\rbrack.
\label{eq:rel2}
\end{eqnarray*}
\end{enumerate}
We call $\boxtimes$ the {\it tetrahedron algebra}.
\end{definition}

\noindent
We comment on how $\boxtimes $ is related to
 $\mathfrak{sl}_2$. Recall that 
$\mathfrak{sl}_2$
is the Lie algebra over $\K$
with a basis $h, e^{\pm}$ 
 and Lie bracket
\begin{eqnarray*}
\lbrack h,e^{\pm}\rbrack = \pm 2e^{\pm},
\qquad
\qquad
\lbrack e^+,e^-\rbrack = h.
\end{eqnarray*}
Define
\begin{eqnarray*}
x=h,
\qquad \qquad 
y=2e^+-h,
\qquad \qquad 
z=-2e^- -h.
\end{eqnarray*}
Then $x,y,z$ is a basis for 
$\mathfrak{sl}_2$ and 
\begin{eqnarray}
\label{eq:equit}
\lbrack
x,y
\rbrack =  2x+2y,
\qquad 
\qquad
\lbrack
y,z
\rbrack =  2y+2z,
\qquad 
\qquad
\lbrack
z,x
\rbrack =  2z+2x.
\end{eqnarray}
We call
$x,y,z$ the {\it equitable basis} for 
$\mathfrak{sl}_2$.

\begin{proposition} {\rm \cite[Corollary 12.1]{HT}}
\label{lem:sl2inj}
Let $h,i,j$ denote mutually distinct elements of $\I$.
Then there exists an injection of Lie algebras from 
$\mathfrak{sl}_2$ to $\boxtimes$ that sends
\begin{eqnarray*}
x \;\rightarrow \; x_{hi},
\qquad \qquad 
y \;\rightarrow  \;x_{ij},
\qquad \qquad 
z \;\rightarrow  \;x_{jh}.
\end{eqnarray*}
\end{proposition}

%\begin{note}
%\rm
%By Proposition 
%\ref{lem:sl2inj}, for $h \in \I$ the subspace
%\begin{eqnarray*}
%\boxtimes^{(h)}:=\mbox{Span}(X_{ij} \;|\;i\not=j,\;i\not=h,\;j\not=h)
%\end{eqnarray*}
%is a subalgebra of $\boxtimes$ that is isomorphic to $\mathfrak{sl}_2$.
%Moreover $\boxtimes$ is generated by the $\boxtimes^{(h)}$, $h \in \I$.
%\end{note}

\section{The ${\mathfrak{sl}}_2$ loop algebra}
In this section we discuss how $\boxtimes$ is related to
the 
 ${\mathfrak{sl}}_2$ loop algebra. We begin with a definition.

\begin{definition}
\label{lem:loopkac}
{\rm \cite[p.~100]{kac}} \rm
The loop algebra 
$L(\mathfrak{sl}_2)$ is the Lie algebra over
$\K$ that has generators $h_i, e^{\pm}_i,$ $i\in \lbrace 0,1\rbrace$
and the following relations:
\begin{eqnarray*}
h_0+h_1 &=& 0,
\\
\lbrack h_i, e_i^{\pm} \rbrack &=& \pm 2 e^{\pm}_i,
\\
\lbrack h_i, e^{\pm}_j \rbrack &=& \mp 2 e^{\pm}_j, \quad \qquad i \not=j,
\\
\lbrack e^{+}_i, e^{-}_i \rbrack &=&h_i,
\\
\lbrack e^{\pm}_0, e^{\mp}_1 \rbrack &=&0,
\\
\lbrack e^{\pm}_i, \lbrack e^{\pm}_i, \lbrack e^{\pm}_i, e^{\pm}_j \rbrack \rbrack\rbrack
&=& 0, \qquad \qquad \quad i \not=j.
\end{eqnarray*}
\end{definition}

%%%%%%%save below%%%%%%%%%%%%
%An isomorphism is given by
%\begin{eqnarray*}
%&&e_1 \;\rightarrow \;e\otimes 1,
%\qquad \qquad  \;
%f_1 \;\rightarrow \; f\otimes 1,
%\qquad \qquad \;\; \;\;
%h_1 \;\rightarrow \; h\otimes 1,
%\\
%&&e_0 \;\rightarrow \; f \otimes T,
%\qquad \qquad 
%f_0 \;\rightarrow \; e\otimes T^{-1}, 
%\qquad \qquad 
%h_0 \;\rightarrow \;-h \otimes 1.
%\end{eqnarray*}
%%%%%%%%%%%save above%%%%%%%%%%

\noindent The following  presentation of $L(\mathfrak{sl}_2)$ will
be useful.

\begin{lemma}
\label{lem:loop2}
{\rm \cite[Lemma 5.3]{HT}}
$L(\mathfrak{sl}_2)$ is isomorphic to the Lie algebra over
$\K$ that has generators 
$x_i, y_i, z_i,$ $i \in \lbrace 0,1\rbrace$ and the following
relations:
\begin{eqnarray*}
x_0+x_1&=&0,
\\
\lbrack x_i,y_i \rbrack &=& 2x_i+2y_i,
\\
\lbrack y_i,z_i \rbrack &=& 2y_i+2z_i,
\\
\lbrack z_i,x_i \rbrack &=& 2z_i+2x_i,
\\
\lbrack z_i, y_j \rbrack &=& 2z_i+2y_j, \;\; \quad \qquad i\not=j,
\\
\lbrack y_i,
\lbrack y_i,
\lbrack y_i,
y_j\rbrack \rbrack \rbrack &=& 
4 \lbrack y_i,
y_j\rbrack,      \qquad \qquad i \not=j,
\label{eq:ydg1}
\\
\lbrack z_i,
\lbrack z_i,
\lbrack z_i,
z_j\rbrack \rbrack \rbrack &=& 
4 \lbrack z_i,
z_j\rbrack, \qquad \qquad i \not=j.
\label{eq:ydg2}
\end{eqnarray*}
An isomorphism with the presentation in Lemma
\ref{lem:loopkac} is given by
\begin{eqnarray*}
x_i \;\rightarrow \; h_i,\qquad \qquad 
y_i \;\rightarrow \;2e^{+}_i-h_i,
\qquad \qquad  
z_i \;\rightarrow \;-2e^{-}_i-h_i.
\end{eqnarray*}
The inverse of this isomorphism is given by
\begin{eqnarray*}
h_i \;\rightarrow \;x_i, \qquad \qquad 
e^{+}_i \;\rightarrow \;(x_i+y_i)/2,
\qquad \qquad 
e^{-}_i \;\rightarrow \;-(z_i+x_i)/2.
\end{eqnarray*}
\end{lemma}
\noindent {\it Proof:} 
One routinely checks that each map is a homomorphism
of Lie algebras and that the maps are inverses.
It follows that each map is an isomorphism of Lie algebras.
\hfill $\Box $ \\

\begin{proposition}
{\rm \cite[Corollary 12.3]{HT}}
\label{prop:loopinj}
Let $h,i,j,k$ denote mutually distinct  elements of $\I$.
Then there exists an injection of Lie algebras
from $L(\mathfrak{sl}_2)$ to $\boxtimes$ that sends
\begin{eqnarray*}
&&x_1\;\rightarrow x_{hi},
\qquad \qquad 
y_1\;\rightarrow x_{ij},
\qquad \qquad 
z_1\;\rightarrow x_{jh},
\\
&&x_0\;\rightarrow x_{ih},
\qquad \qquad 
y_0\;\rightarrow x_{hk},
\qquad \qquad 
z_0\;\rightarrow x_{ki}.
\end{eqnarray*}
\end{proposition}

\section{The Onsager algebra}

\noindent In this section we discuss how
$\boxtimes$ is related to the Onsager
algebra.

\begin{definition}\rm
\label{lem:onsalt}
{\rm \cite{Onsager}, \cite{perk}}
The {\it Onsager algebra $\mathcal O$} is
the Lie algebra over
$\K$ that has generators $x,y$ and relations
\begin{eqnarray}
\lbrack x,
\lbrack x,
\lbrack x,
y\rbrack \rbrack \rbrack &=& 
4 \lbrack x,
y\rbrack,
\label{eq:dg1}
\\
\lbrack y,
\lbrack y,
\lbrack y,
x\rbrack \rbrack \rbrack &=& 
4 \lbrack y,
x\rbrack.
\label{eq:dg2}
\end{eqnarray}
\end{definition}

\begin{definition}
\rm
Referring to Definition 
\ref{lem:onsalt},
we call $x,y$ the {\it standard generators}
for ${\mathcal O}$.
\end{definition}

\noindent We refer the reader to 
\cite{DateRoan2},
\cite{Davfirst},
\cite{Da},
%\cite{vg3},
\cite{perk},
\cite{roan}
for a mathematical treatment of 
the Onsager algebra.
For connections to 
solvable lattice models see
\cite{CKOnsn},
\cite{Albert},
\cite{agmpy},
\cite{auyangandperk},
%\cite{auyangperk2},
\cite{auyang2},
\cite{Baz},
\cite{McCoyd},
\cite{Dolgra},
\cite{vongehlen},
\cite{Klish1},
\cite{Klish2},
\cite{Klish3},
\cite{Lee},
\cite{uglov1}.
%tridiagonal pairs
%\cite{hartwig},
%  \cite{TD00},
%  \cite{LS99},
%   \cite{qSerre},
%   \cite{aw}
%   and partially orthogonal polynomials
%   \cite{vg1},
%   \cite{vg2}.

\begin{proposition}{\rm \cite[Corollary 12.2]{HT}}
\label{prop:onshom}
Let $h,i,j,k$ denote mutually distinct elements of $\I$.
Then there exists an injection of Lie algebras from
$\mathcal O$ to $\boxtimes$ that sends
\begin{eqnarray*}
x \rightarrow x_{hi}, 
\qquad \qquad y \rightarrow x_{jk}.
\end{eqnarray*}
\end{proposition}

\noindent We have a remark.
Let $\Omega $ (resp. $\Omega^{\prime}$)
 (resp. $\Omega^{\prime \prime}$)
denote the subalgebra of $\boxtimes$ generated by
$x_{01}$ and $x_{23}$
(resp.
$x_{02}$ and $x_{31}$)
(resp.
$x_{03}$ and $x_{12}$).
By Proposition
\ref{prop:onshom} each of
 $\Omega$, ${\Omega}^{\prime}$,
  ${\Omega}^{\prime \prime}$ is
isomorphic to $\mathcal O$. By \cite[Theorem 11.6]{HT} we have
\begin{eqnarray*}
\boxtimes =  {\Omega} +{\Omega}^{\prime}+
{\Omega}^{\prime \prime} \qquad \qquad
(\mbox{direct sum}).
\end{eqnarray*}

\section{The finite dimensional irreducible $\boxtimes$-modules}

\noindent 
Throughout this section  assume that the field $\K$ is 
algebraically closed with characteristic 0.

\medskip
\noindent In \cite{Ha}
Hartwig classifies the finite dimensional irreducible
$\boxtimes$-modules. He does this by relating them
to the finite dimensional irreducible ${\mathcal O}$-modules
which were previously classified by Davies \cite{Da}; see
also 
\cite{DateRoan2},
\cite{roan}.
We will state Hartwig's main results after a few comments.

\medskip
\noindent Let $V$ denote a finite dimensional irreducible
${\mathcal O}$-module. Define a scalar $\alpha$ (resp. $\alpha^*$)
to be  $(\mbox{\rm dim}(V))^{-1}$ times
the trace of $x$ (resp. $y$) on $V$, where
$x,y$ are the standard generators for ${\mathcal O}$.
We call the ordered pair $(\alpha, \alpha^*)$ the
{\it type} of $V$. Replacing $x$ and $y$ by $x-\alpha I$ and 
$y-\alpha^*I$ respectively, the type becomes $(0,0)$.

\begin{theorem}
\label{thm:Hart1}
{\rm (Hartwig \cite[Theorem  1.2]{Ha})}
Let $V$ denote a finite dimensional irreducible
$\boxtimes$-module. Then there exists 
a unique ${\mathcal O}$-module structure on $V$
such that the standard generators $x$ and $y$ act as
$x_{01}$ and $x_{23}$ respectively.
This ${\mathcal O}$-module structure is irreducible
of type $(0,0)$.
\end{theorem}

\begin{theorem}
\label{thm:Hart2}
{\rm (Hartwig \cite[Theorem 1.3]{Ha})}
Let $V$ denote a finite dimensional irreducible
${\mathcal O}$-module of type $(0,0)$. Then
there exists a unique $\boxtimes$-module structure
on $V$ such that the standard generators $x$ and $y$
act as $x_{01}$  and $x_{23}$ respectively.
This $\boxtimes$-module structure is irreducible.
\end{theorem}

\begin{remark}
{\rm (Hartwig \cite[Remark 1.4]{Ha})}
\rm
Combining Theorem 
\ref{thm:Hart1} and
 Theorem 
\ref{thm:Hart2} we obtain a bijection between
the following two sets:
\begin{enumerate}
\item[{\rm (i)}]  
the isomorphism classes of finite dimensional irreducible
$\boxtimes$-modules.
\item[{\rm (ii)}] 
the isomorphism classes of finite dimensional irreducible
${\mathcal O}$-modules of type $(0,0)$.
\end{enumerate}
\end{remark}

\section{The $q$-tetrahedron algebra}

\noindent We are now ready to introduce the 
$q$-tetrahedron algebra. 
Until further notice  assume that the  field $\K$ is arbitrary.
We fix a nonzero scalar $q \in \K$ such that $q^2\not=1$ and define
\begin{eqnarray}
\lbrack n \rbrack_q = \frac{q^n-q^{-n}}{q-q^{-1}}
\qquad \qquad n = 0,1,2,\ldots 
\label{eq:nbrack}
\end{eqnarray}
We let $\Z_4 = \Z/4\Z$ denote the cyclic group of order 4.

\begin{definition}
\label{def:qtet}
\rm
Let $\boxtimes_q$ denote the unital associative $\K$-algebra that has
generators 
\begin{eqnarray*}
\lbrace x_{ij}\;|\; i,j \in \Z_4,\;j-i=1 \;\mbox{or} \;j-i=2\rbrace
\end{eqnarray*}
and the following relations:
\begin{enumerate}
\item For $i,j\in \Z_4$ such that $j-i=2$,
\begin{eqnarray*}
x_{ij}x_{ji} = 1.
\label{eq:qrel0}
\end{eqnarray*}
\item For $h,i,j\in \Z_4$ such that the pair $(i-h,j-i)$ is one of
$(1,1), (1,2), (2,1)$,
\begin{eqnarray*}
\frac{qx_{hi}x_{ij}-q^{-1}x_{ij}x_{hi}}{q-q^{-1}}=1.
\label{eq:qrel1}
\end{eqnarray*}
\item For $h,i,j,k\in \Z_4$ such that $i-h=j-i=k-j=1$,
\begin{eqnarray}
\label{eq:qserre}
x_{hi}^3x_{jk} -
\lbrack 3 \rbrack_q
x_{hi}^2x_{jk}x_{hi} +
\lbrack 3 \rbrack_q
x_{hi}x_{jk}x_{hi}^2- 
x_{jk}x_{hi}^3=0. 
\end{eqnarray}
\end{enumerate}
We call $\boxtimes_q$ the {\it $q$-tetrahedron algebra}.
\end{definition}
\begin{note}\rm
The equations (\ref{eq:qserre}) are the cubic $q$-Serre relations.
\end{note}
\noindent
We make some observations.

\begin{lemma}
\label{lem:rho}
There exists a $\K$-algebra automorphism $\rho$ of $\boxtimes_q$
that sends each generator $x_{ij}$ to $x_{i+1,j+1}$.
Moreover 
 $\rho^4=1$.
\end{lemma}

\begin{lemma}
\label{lem:omega}
There exists a $\K$-algebra antiautomorphism $\omega$ of 
 $\boxtimes_q$
that sends
\begin{eqnarray*}
&&x_{01} \to x_{01},\quad 
x_{12} \to x_{30}, \quad x_{23} \to x_{23}, \quad x_{30} \to x_{12},
\\
&&x_{02} \to x_{31},\quad 
x_{13} \to x_{20}, \quad x_{20} \to x_{13}, \quad x_{31} \to x_{02}.
\end{eqnarray*}
Moreover $\omega^2=1$.
\end{lemma}

\begin{lemma}
\label{lem:flip}
There exists a $\K$-algebra automorphism of 
 $\boxtimes_q$
that sends each generator $x_{ij}$ to $-x_{ij}$.
%This automorphism has order 2.
\end{lemma}

\section{The algebra $U_q(\mathfrak{sl}_2)$}

\noindent In this section we discuss how the algebra
$\boxtimes_q$ is related to 
$U_q(\mathfrak{sl}_2)$.
We start with a definition.

\begin{definition} 
\label{def:uq}
\rm
\cite[p.~122]{Kassel}
Let $U_q(\mathfrak{sl}_2)$
denote the unital associative $\K$-algebra 
with
generators $K^{\pm 1}$, $e^{\pm}$
and the following relations:
\begin{eqnarray*}
KK^{-1} &=& 
K^{-1}K =  1,
\label{eq:buq1}
\\
Ke^{\pm}K^{-1} &=& q^{\pm 2}e^{\pm},
\label{eq:buq2}
\\
\lbrack e^+,e^-\rbrack  &=& \frac{K-K^{-1}}{q-q^{-1}}.
\label{eq:buq4}
\end{eqnarray*}
\end{definition}

%\noindent 
%We call $k^{\pm 1}, e,f$
%the {\it Chevalley generators} for
%$U_q(\mathfrak{sl}_2)$.

\noindent The following presentation of
$U_q(\mathfrak{sl}_2)$ will be useful.
%given
%in 
%Definition
%\ref{def:uq}
%the generators $k^{\pm 1}$
%and the generators $e,f$ play a very different
%role. 
%We now recall a presentation for 
%$U_q(\mathfrak{sl}_2)$ whose
%generators
%are on a more equal footing.  

\begin{lemma}
\label{thm:uq2}
{\rm \cite[Theorem 2.1]{equit1}}
The algebra
$U_q(\mathfrak{sl}_2)$ 
 is isomorphic to
the unital associative $\K$-algebra 
with
generators 
$x^{\pm 1}$, $y$, $z$
and the following relations:
\begin{eqnarray*}
xx^{-1} = 
x^{-1}x &=&  1,
\label{eq:2buq1}
\\
\frac{qxy-q^{-1}yx}{q-q^{-1}}&=&1,
\label{eq:2buq2}
\\
\frac{qyz-q^{-1}zy}{q-q^{-1}}&=&1,
\label{eq:2buq3}
\\
\frac{qzx-q^{-1}xz}{q-q^{-1}}&=&1.
\label{eq:2buq4}
\end{eqnarray*}
An isomorphism with the presentation in Definition
\ref{def:uq} is given by:
\begin{eqnarray*}
\label{eq:iso1}
x^{{\pm}1} &\rightarrow & K^{{\pm}1},\\
\label{eq:iso2}
%y &\rightarrow & k^{-1}+f(q-q^{-1}), \\
%\label{eq:iso3}
%z &\rightarrow & k^{-1}-k^{-1}eq(q-q^{-1}).
y &\rightarrow & K^{-1}+e^-, \\
\label{eq:iso3}
z &\rightarrow & K^{-1}-K^{-1}e^+q(q-q^{-1})^2.
\end{eqnarray*}
The inverse of this isomorphism is given by:
\begin{eqnarray*}
\label{eq:iso1inv}
K^{{\pm}1} &\rightarrow & x^{{\pm}1},\\
\label{eq:iso2inv}
%f &\rightarrow & (y-x^{-1})(q-q^{-1})^{-1}, \\
%\label{eq:iso3inv}
%e &\rightarrow & (1-xz)q^{-1}(q-q^{-1})^{-1}.
e^- &\rightarrow & y-x^{-1}, \\
\label{eq:iso3inv}
e^+ &\rightarrow & (1-xz)q^{-1}(q-q^{-1})^{-2}.
\end{eqnarray*}
\end{lemma}
\noindent {\it Proof:} One readily checks that
each map is a homomorphism of $\K$-algebras and that
the maps are inverses.
It follows that each map is an isomorphism of $\K$-algebras.
\hfill $\Box $ \\

\begin{definition}
\rm
\cite{equit1}
 Referring to Lemma  
\ref{thm:uq2},
We call 
$x^{\pm 1}, y,z$ the
{\it equitable generators} for
$U_q(\mathfrak{sl}_2)$. 
\end{definition}

\begin{proposition}
\label{lem:uqinj}
For $i \in \Z_4$ there exists a $\K$-algebra homomorphism
from
$U_q(\mathfrak{sl}_2)$ to
 $\boxtimes_q$ that sends 
\begin{eqnarray*}
x\to x_{i,i+2},
\quad 
x^{-1}\to x_{i+2,i},\quad
y\to x_{i+2,i+3},
\quad 
z\to x_{i+3,i}.
\end{eqnarray*}
\end{proposition}
\noindent {\it Proof:}
Compare the defining relations for 
$U_q(\mathfrak{sl}_2)$
given in Lemma  
\ref{thm:uq2}
with the relations 
in Definition
\ref{def:qtet}(i),(ii).
\hfill $\Box $ 

\begin{conjecture}
\rm
The map in Proposition 
\ref{lem:uqinj} is an injection.
\end{conjecture}

%\begin{note}
%\rm
%By Lemma \ref{lem:uqinj}
%for $h \in \Z_4$ the subalgebra
%$\boxtimes_q^{(h)}$ generated by
%\begin{eqnarray*}
%x_{h+1,h+2},\quad x_{h+2,h+3},\quad
%x_{h+3,h+1},\quad x_{h+1,h+3}
%\end{eqnarray*}
%is a homomorphic image of 
%$U_q(\mathfrak{sl}_2)$. By Definition
%\ref{def:qtet} the algebra $\boxtimes_q$ is generated by the
%$\boxtimes_q^{(h)}$, $h \in \Z_4$.
%\end{note}

\section{
The 
$U_q(\mathfrak{sl}_2)$
loop algebra 
}

\noindent In this section we consider how
$\boxtimes_q$ is related to the
$U_q(\mathfrak{sl}_2)$ loop algebra.
We start with a definition.

\begin{definition} 
\label{def:loopuq}
\rm
\cite[p.~266]{charp} 
Let  
$U_q(L(\mathfrak{sl}_2))$ denote the unital associative $\K$-algebra 
with
generators $K_i$,
 $e^{\pm}_i$,
$i\in \lbrace 0,1\rbrace $
and the following relations:
\begin{eqnarray*}
K_0K_1&=& K_1K_0=1,
\label{eq:lbuq2}
\\
K_ie^{\pm}_iK^{-1}_i &=& q^{{\pm}2}e^{\pm}_i,
\label{eq:lbuq3}
\\
K_ie^{\pm}_jK^{-1}_i &=& q^{{\mp}2}e^{\pm}_j, \qquad i\not=j,
\label{eq:lbuq4}
\\
\lbrack e^+_i, e^-_i\rbrack &=& {{K_i-K^{-1}_i}\over {q-q^{-1}}},
\label{eq:lbuq5}
\\
\lbrack e^{\pm}_0, e^{\mp}_1\rbrack &=& 0,
\label{eq:lbuq6}
\end{eqnarray*}
\begin{eqnarray*}
(e^{\pm}_i)^3e^{\pm}_j -  
\lbrack 3 \rbrack_q (e^{\pm}_i)^2e^{\pm}_j e^{\pm}_i 
+\lbrack 3 \rbrack_q e^{\pm}_ie^{\pm}_j (e^{\pm}_i)^2 - 
e^{\pm}_j (e^{\pm}_i)^3 =0, \qquad i\not=j.
\label{eq:lbuq7}
\end{eqnarray*}
We call $U_q(L({\mathfrak{sl}_2}))$ the
  $U_q(\mathfrak{sl}_2)$
{\it loop algebra}.
%We call $e^{\pm}_i$, $K_i^{{\pm}1}$, $i\in \lbrace 0,1\rbrace $
%the {\it Chevalley generators} for
%$U_q({\widehat{sl}}_2)$.
\end{definition}

\noindent The following presentation of
$U_q(L(\mathfrak{sl}_2))$ will be useful.

\begin{theorem}
\label{thm:luq2} 
{\rm (\cite[Theorem 2.1]{tdanduq},
\cite{equit2})}
The loop algebra
$U_q(L(\mathfrak{ sl}_2))$ is isomorphic to
the unital associative $\K$-algebra 
with
generators $x_i,y_i,z_i$, $i\in \lbrace 0,1\rbrace $
and the following relations:
\begin{eqnarray*}
x_0x_1=x_1x_0&=&1, 
\label{eq:l2buq2}
\\
\frac{q x_iy_i-q^{-1}y_ix_i}{q-q^{-1}} &=& 1,
\label{eq:l2buq3}
\\
\frac{q y_iz_i-q^{-1}z_iy_i}{q-q^{-1}} &=& 1,
\label{eq:l2buq4}
\\
\frac{q z_ix_i-q^{-1}x_iz_i}{q-q^{-1}} &=& 1,
\label{eq:l2buq5}
\\
\frac{q z_iy_j-q^{-1}y_jz_i}{q-q^{-1}} &=&1,
\qquad i\not=j,
\label{eq:l2buq6}
\end{eqnarray*}
\begin{eqnarray*}
y_i^3y_j -  
\lbrack 3 \rbrack_q y_i^2y_j y_i 
+\lbrack 3 \rbrack_q y_iy_j y_i^2 - 
y_j y_i^3 =0, \qquad i\not=j,
\label{eq:l2buq7}
\\
z_i^3z_j -  
\lbrack 3 \rbrack_q z_i^2z_j z_i 
+\lbrack 3 \rbrack_q z_iz_j z_i^2 - 
z_j z_i^3 =0, \qquad i\not=j.
\label{eq:l2buq8*}
\end{eqnarray*}
An isomorphism with the presentation in Definition
\ref{def:loopuq} is given by:
\begin{eqnarray*}
\label{eq:liso1}
x_i &\rightarrow & K_i,\\
\label{eq:liso2}
y_i &\rightarrow & K^{-1}_i+e^{-}_i, \\
\label{eq:liso3}
z_i &\rightarrow & K^{-1}_i-K^{-1}_ie^{+}_iq(q-q^{-1})^2.
\end{eqnarray*}
The inverse of this isomorphism is given by:
\begin{eqnarray*}
\label{eq:liso1inv}
K_i &\rightarrow & x_i,\\
\label{eq:liso2inv}
e^-_i &\rightarrow & y_i-x^{-1}_i, \\
\label{eq:liso3inv}
e^+_i &\rightarrow & (1-x_i z_i)q^{-1}(q-q^{-1})^{-2}.
\end{eqnarray*}
\end{theorem}
\noindent {\it Proof:} One readily checks that
each map is a homomorphism of $\K$-algebras and that
the maps are inverses.
It follows that each map is an isomorphism of $\K$-algebras.
\hfill $\Box $ \\

%\begin{definition}
%\label{def:alt}
%\rm 
%With reference to Theorem
%\ref{thm:uq2}
%we call
% $y^{\pm}_i$, $k_i$, $i\in \lbrace 0,1\rbrace $
%the {\it equitable generators} of
%$L(U_q({\mathfrak{sl}}_2))$.
%\end{definition}

\begin{proposition}
\label{prop:loop}
For $i \in \Z_4$ there exists a $\K$-algebra homomorphism
from 
$U_q(L({\mathfrak{sl}}_2))$ to $\boxtimes_q$
that sends
\begin{eqnarray*}
&&
x_1 \to x_{i,i+2},\quad 
y_1\to x_{i+2,i+3},\quad
z_1 \to x_{i+3,i},
\\
&&x_0 \to x_{i+2,i}, \quad 
y_0 \to x_{i,i+1},\quad
z_0 \to x_{i+1,i+2}.
\end{eqnarray*}
\end{proposition}
\noindent {\it Proof:} 
Compare the defining relations for
$U_q(L({\mathfrak{sl}}_2))$ given in
Theorem
\ref{thm:luq2} with the relations
in Definition
\ref{def:qtet}.
\hfill $\Box $

\begin{conjecture}
\rm
The map in Proposition
\ref{prop:loop} is an injection.
\end{conjecture}

\section{The algebra ${\mathcal A}_q$}

In this section we discuss how $\boxtimes_q$ is 
related to the algebra ${\mathcal A}_q$.

\begin{definition}
\label{defa}
\rm
Let ${\mathcal A}_q$ denote the unital associative $\K$-algebra
defined by generators $x,y$
and relations
\begin{eqnarray}
\label{eq:qs1}
x^3y-\lbrack 3\rbrack_q x^2yx
+\lbrack 3\rbrack_q xyx^2
-yx^3&=&0,
\\
\label{eq:qs2}
y^3x-\lbrack 3\rbrack_q y^2xy
+\lbrack 3\rbrack_q yxy^2
-xy^3&=&0.
\end{eqnarray}
%We call $x,y$ the {\it standard generators} for ${\mathcal A}_q$.
\end{definition}

\begin{definition} \rm
Referring to Definition
\ref{defa}, we call $x,y$ the {\it standard generators}
for 
${\mathcal A}_q$.
\end{definition}

\begin{note}\rm
{\rm \cite[Corollary 3.2.6]{lusztig}}
The algebra ${\mathcal A}_q$ is often called the
{\it positive part of 
$U_q(\widehat{ \mathfrak{sl}}_2)$.}
\end{note}

\begin{proposition}
\label{lem:aqinj}
For $i \in \Z_4$ there exists a homomorphism of $\K$-algebras
from ${\mathcal A}_q$ to $\boxtimes_q$
that sends the standard generators
$x, y$ to $x_{i,i+1}, x_{i+2,i+3}$ respectively.
\end{proposition}
\noindent {\it Proof:}
Compare the relations
(\ref{eq:qs1}),
(\ref{eq:qs2})
with the relations
(\ref{eq:qserre}).
\hfill $\Box $ \\

\begin{conjecture}
The map in Proposition 
\ref{lem:aqinj} is an injection.
\end{conjecture}

\section{The finite dimensional
irreducible $\boxtimes_q$-modules}

\noindent For the rest of this paper assume
that $q$ is not a root of 1, and
the field $\K$ is algebraically closed.

\medskip
\noindent Our next goal is to describe  
 the finite dimensional irreducible
$\boxtimes_q$-modules. 
As part of this description we relate these modules
to a type of ${\mathcal A}_q$-module that
we  discussed in \cite{NN}.
In order to define this type of 
${\mathcal A}_q$-module 
we recall a concept. Let $V$ denote a 
finite-dimensional vector space over $\K$. 
A linear transformation $A:V\to V$ is said
to be {\it nilpotent} whenever
there exists a positive integer $n$ such that
$A^n=0$. 

\begin{definition}
\label{def:nn}
\rm
\cite[Definition 1.3]{NN}
Let $V$ denote a finite-dimensional
${\mathcal A}_q$-module.
We say this module is {\it NonNil}
whenever the standard generators $x, y$ are not nilpotent
on $V$.
\end{definition}

\begin{note}
\rm In \cite[Theorems 1.6, 1.7]{NN} we 
classified up to isomorphism the 
NonNil finite-dimensional
irreducible ${\mathcal A}_q$-modules 
assuming the field $\K$ has characteristic $0$ in addition
to being algebraically closed.
\end{note}

\noindent We now describe how the
finite dimensional
irreducible $\boxtimes_q$-modules are
related to the NonNil finite dimensional
irreduble ${\mathcal A}_q$-modules.
We begin with a few comments.
Let $V$ denote a
NonNil
finite-dimensional
irreducible
${\mathcal A}_q$-module.
By \cite[Corollary 2.8]{NN} 
 the standard generators $x,y$ are
semisimple on
$V$. Moreover
there exist an integer $d\geq 0$ and nonzero scalars $\alpha, \alpha^* \in \K$
such that the set of distinct eigenvalues of $x$ (resp. $y$) on $V$
is $\lbrace \alpha q^{d-2n} \,|\, 0 \leq n \leq d\rbrace$
(resp. $\lbrace \alpha^*q^{d-2n} \,|\, 0 \leq n \leq d\rbrace$).
We call the ordered pair $(\alpha,\alpha^*)$ the {\it type}
of $V$.
Replacing $x,y$ by $x/\alpha, y/\alpha^*$ the type becomes $(1,1)$.
Now let $V$ denote a finite dimensional irreducible 
$\boxtimes_q$-module. As we will see, there exist an
integer $d\geq 0$ and a scalar $\varepsilon \in \lbrace 1,-1\rbrace$
such that for each generator $x_{ij}$ the action on
$V$ is semisimple with eigenvalues $\lbrace \varepsilon q^{d-2n}\;|\;0\leq n\leq d\rbrace$. We call $\varepsilon $ the {\it type} of $V$.
Replacing each generator $x_{ij}$ by $\varepsilon x_{ij}$ the
type becomes 1. The main results of the present paper are contained
in the following two theorems and subsequent remark. 

\begin{theorem}
\label{thm:2}
Let $V$ denote a 
finite-dimensional
irreducible 
$\boxtimes_q$-module of type 1.
Then there exists
a unique
${\mathcal A}_q$-module 
structure on 
$V$ such that the standard generators
$x$ and  $y$ act as $x_{01}$  
and $x_{23}$ respectively.
This 
${\mathcal A}_q$-module 
is NonNil irreducible 
of type  $(1,1)$.
\end{theorem}

\begin{theorem}
\label{thm:1}
Let $V$ denote a
NonNil 
finite-dimensional
irreducible
${\mathcal A}_q$-module of type $(1,1)$.
Then there exists
a unique
$\boxtimes_q$-module
structure on 
$V$ such that the standard generators
$x$ and  $y$ act as $x_{01}$  
and $x_{23}$ respectively.
This 
$\boxtimes_q$-module structure
is irreducible and type $1$.
\end{theorem}

\begin{remark}
\rm Combining Theorem
\ref{thm:2} and
Theorem
\ref{thm:1}
we obtain a bijection between the following two sets:
\begin{enumerate}
\item the isomorphism classes of 
finite-dimensional 
irreducible 
$\boxtimes_q$-modules
of type $1$;
\item the isomorphism classes of NonNil
finite-dimensional 
irreducible 
${\mathcal A}_q$-modules of type $(1,1)$.
\end{enumerate}
\end{remark}

\noindent 
We prove Theorems 
\ref{thm:2} and
\ref{thm:1} in
Sections 17 and 18, respectively.
In Sections 11--16 we will obtain some 
results used in these proofs.

\section{Some linear algebra}

\noindent In this section we obtain some
linear algebraic results
that will help us describe the finite dimensional irreducible
$\boxtimes_q$-modules.

\medskip
\noindent 
We will use the following concepts.
Let $V$ denote a finite-dimensional
vector space over $\K$ 
and let $A:V\to V$ denote a linear transformation.
For $\theta \in \K$ we define
\begin{eqnarray*}
V_A(\theta) = \lbrace v \in V \,|\,Av = \theta v\rbrace.
\end{eqnarray*}
Observe that $\theta$ is an eigenvalue of $A$ if and only
if $V_A(\theta)\not=0$, and in this case
$V_A(\theta)$ is the corresponding eigenspace.
The sum
$\sum_{\theta \in \K} V_A(\theta)$ is direct.
Moreover 
this sum is equal to $V$ 
if and only if
$A$ is semisimple.

\begin{lemma}
\label{lem:key}
Let $V$ denote a finite-dimensional
vector space over $\K$.
Let $A:V\to V$ and $B:V\to V$ denote linear transformations.
Then for all nonzero $\theta \in \K$ the following 
are equivalent:
\begin{enumerate}
\item[{\rm (i)}]  
The expression
$
A^3B-\lbrack 3\rbrack_q A^2BA
+\lbrack 3\rbrack_q ABA^2
-BA^3
$
vanishes on $V_A(\theta)$.
\item[{\rm (ii)}]
$B V_A(\theta) \subseteq 
 V_A(q^2\theta)
+
 V_A(\theta)
+
V_A(q^{-2}\theta)
$.
\end{enumerate}
\end{lemma}
\noindent {\it Proof:} 
For $v \in V_A(\theta)$ we
have 
\begin{eqnarray*}
&&(A^3B-\lbrack 3\rbrack_q A^2BA
+\lbrack 3\rbrack_q ABA^2
-BA^3)v
\\
&& \qquad = 
(A^3-\theta \lbrack 3\rbrack_q A^2
+\theta^2\lbrack 3\rbrack_q A
-\theta^3I)Bv
\qquad \qquad \mbox{Since}\;Av=\theta v
\\
&&\qquad =
(A-q^2 \theta I)
(A- \theta I)
(A-q^{-2} \theta I)Bv,
\end{eqnarray*}
where $I:V\to V$ is the identity map.
The scalars $q^{2}\theta, \theta, q^{-2}\theta$
are mutually distinct since $\theta\not=0$ and
since $q$ is not a root of $1$.
The result follows.
\hfill $\Box $

\begin{lemma}
\label{lem:qweyl}
Let $V$ denote a vector space over $\K$ with finite positive
dimension. Let $A:V\to V$ and $B:V\to V$ denote linear transformations.
Then for all nonzero $\theta \in \K$ the following are equivalent:
\begin{enumerate}
\item[{\rm (i)}]  
The expression
$
qAB-q^{-1}BA-(q-q^{-1})I
$
vanishes on $V_A(\theta)$.
\item[{\rm (ii)}]
$(B -\theta^{-1}I)V_A(\theta) \subseteq 
 V_A(q^{-2}\theta)
$.
\end{enumerate}
\end{lemma}
\noindent {\it Proof:} 
For $v \in V_A(\theta)$ we
have 
\begin{eqnarray*}
(qAB-q^{-1}BA-(q-q^{-1})I)v=
q(A-q^{-2}\theta I)(B-\theta^{-1}I)v
\end{eqnarray*}
and the result follows.
\hfill $\Box $ \\ 

\noindent For later use we give a second version of Lemma
\ref{lem:qweyl}.

\begin{lemma}
\label{lem:qweyl2}
Let $V$ denote a vector space over $\K$ with finite positive
dimension. Let $A:V\to V$ and $B:V\to V$ denote linear transformations.
Then for all nonzero $\theta \in \K$ the following are equivalent:
\begin{enumerate}
\item[{\rm (i)}]  
The expression
$
qAB-q^{-1}BA-(q-q^{-1})I
$
vanishes on $V_B(\theta)$.
\item[{\rm (ii)}]
$(A -\theta^{-1}I)V_B(\theta) \subseteq 
 V_B(q^{2}\theta)
$.
\end{enumerate}
\end{lemma}
\noindent {\it Proof:} 
In Lemma
\ref{lem:qweyl} replace
$(A,B,q)$ by $(B,A,q^{-1})$.
\hfill $\Box $  \\

\begin{lemma}
\label{lem:qweyl3}
Let $V$ denote a vector space over $\K$ with finite positive
dimension. Let $A:V\to V$ and $B:V\to V$ denote linear transformations
such that
\begin{eqnarray*}
\frac{qAB-q^{-1}BA}{q-q^{-1}}=I.
\end{eqnarray*}
Then for all nonzero $\theta \in \K$,
\begin{eqnarray}
\sum_{n=0}^\infty V_A(q^{-2n}\theta)=
\sum_{n=0}^\infty V_B(q^{2n}\theta^{-1}).
\label{eq:sumsum}
\end{eqnarray}
\end{lemma}
\noindent {\it Proof:} 
Let $S_A$ (resp. $S_B$) denote the sum on the left (resp. right)
in 
(\ref{eq:sumsum}).
We show $S_A=S_B$. Since $V$ has finite dimension
there exists an integer $t\geq 0$ 
such that 
$V_A(q^{-2n}\theta)=0$ and 
$V_B(q^{2n}\theta^{-1})=0$ for $n>t$.
So
$S_A=\sum_{n=0}^t V_A(q^{-2n}\theta)$ and
$S_B=\sum_{n=0}^t V_B(q^{2n}\theta^{-1})$.
By construction $S_B$ is the set of vectors in $V$ on which the
product
\begin{eqnarray}
\prod_{n=0}^t (B-q^{2n}\theta^{-1}I)
\label{eq:bprod}
\end{eqnarray}
is zero. Using Lemma
\ref{lem:qweyl} we find 
(\ref{eq:bprod}) is zero on $S_A$ so $S_A \subseteq S_B$.
By construction $S_A$ is the set of vectors in $V$ on which the
product
\begin{eqnarray}
\prod_{n=0}^t (A-q^{-2n}\theta I)
\label{eq:aprod}
\end{eqnarray}
is zero. Using Lemma
\ref{lem:qweyl2} we find 
(\ref{eq:aprod}) is zero on $S_B$
so $S_B\subseteq S_A$. We conclude
$S_A=S_B$ and the result follows.
\hfill $\Box $  \\

\begin{lemma}
\label{lem:dim}
With the notation and assumptions of Lemma
\ref{lem:qweyl3}, for all nonzero $\theta \in \K$ we have
\begin{eqnarray}
\label{eq:dimsum}
\mbox{\rm dim}V_A(\theta)=
\mbox{\rm dim}V_B(\theta^{-1}).
\end{eqnarray}
\end{lemma}
\noindent {\it Proof:} 
By Lemma \ref{lem:qweyl3} (with $\theta$ replaced by $q^{-2}\theta$),
\begin{eqnarray}
\sum_{n=1}^\infty V_A(q^{-2n}\theta)=
\sum_{n=1}^\infty V_B(q^{2n}\theta^{-1}).
\label{eq:sumsum2}
\end{eqnarray}
Let $S$ denote the sum on either side of
(\ref{eq:sumsum2}). Comparing
(\ref{eq:sumsum})
and
(\ref{eq:sumsum2})
we find
\begin{eqnarray}
V_A(\theta)+S = V_B(\theta^{-1})+S
\label{eq:sdir}
\end{eqnarray}
and
the sum on either side of
(\ref{eq:sdir}) is direct. The result follows.
\hfill $\Box $  \\

\section{The type of a finite dimensional irreducible $\boxtimes_q$-module}

\noindent Let $V$ denote a finite dimensional irreducible
$\boxtimes_q$-module. In this section we show that each
generator $x_{ij}$ of $\boxtimes_q$
is semisimple on $V$. We also find the eigenvalues. 
Using these eigenvalues we associate with $V$
a parameter called the {\it type}.

\medskip

\noindent Before proceeding we refine our notation.
\begin{definition}
\label{def:es}
\rm
Let $V$ denote a finite dimensional irreducible
$\boxtimes_q$-module. For each generator $x_{ij}$ of $\boxtimes_q$
and for each $\theta \in \K$ we write
\begin{eqnarray*}
V_{ij}(\theta)= \lbrace v \in V\;|\;x_{ij}v = \theta v\rbrace.
\end{eqnarray*}
\end{definition}

\begin{lemma}
\label{lem:diminv}
Let $V$ denote a finite dimensional irreducible
$\boxtimes_q$-module and choose a generator $x_{ij}$ of
$\boxtimes_q$.
Then for all nonzero $\theta \in \K$ the spaces
$V_{ij}(\theta)$ and $V_{ij}(\theta^{-1})$ have the same
dimension. This dimension is independent of our choice of generator.
\end{lemma}
\noindent {\it Proof:} 
By Definition
\ref{def:qtet}(ii) we have
\begin{eqnarray}
x_{02}\;\rightarrow \;
x_{23}\;\rightarrow \;
x_{31}\;\rightarrow \;
x_{12}\;\rightarrow \;
x_{20}\;\rightarrow \;
x_{01}\;\rightarrow \;
x_{13}\;\rightarrow \;
x_{30}
\label{eq:chain}
\end{eqnarray}
where $r\rightarrow s$ means
\begin{eqnarray*}
\frac{qrs-q^{-1}sr}{q-q^{-1}}=1.
\end{eqnarray*}
Applying Lemma
\ref{lem:dim}
to each arrow in
(\ref{eq:chain}) we find
\begin{eqnarray}
\label{eq:dimchain}
\mbox{\rm dim}V_{02}(\theta)\;=\;
\mbox{\rm dim}V_{23}(\theta^{-1})\;=\;
\mbox{\rm dim}V_{31}(\theta)\;=\;\cdots \;=\;
\mbox{\rm dim}V_{30}(\theta^{-1})
\end{eqnarray}
and
\begin{eqnarray}
\label{eq:dimchain2}
\mbox{\rm dim}V_{02}(\theta^{-1})\;=\;
\mbox{\rm dim}V_{23}(\theta)\;=\;
\mbox{\rm dim}V_{31}(\theta^{-1})\;=\;\cdots \;=\;
\mbox{\rm dim}V_{30}(\theta).
\end{eqnarray}
Also $V_{02}(\theta)=V_{20}(\theta^{-1})$ since $x_{02}, x_{20}$
are inverses; therefore
\begin{eqnarray}
\label{eq:dimchain3}
\mbox{\rm dim}V_{02}(\theta)&=&
\mbox{\rm dim}V_{20}(\theta^{-1}).
\end{eqnarray}
Combining 
(\ref{eq:dimchain})--(\ref{eq:dimchain3}) we obtain the result.
\hfill $\Box $  \\

\begin{theorem}
\label{thm:diminv}
Let $V$ denote a finite dimensional irreducible
$\boxtimes_q$-module. Then the following hold: 
\begin{enumerate}
\item[{\rm (i)}]  
Each generator $x_{ij}$ of $\boxtimes_q$ is semisimple on $V$.
\item[{\rm (ii)}]
There exist an integer $d\geq 0$ and a scalar $\varepsilon \in \lbrace 1,-1\rbrace$ such that for each generator $x_{ij}$ the set of distinct eigenvalues of
$x_{ij}$ on $V$ is $\lbrace \varepsilon q^{d-2n}\;|\;0 \leq n \leq d\rbrace$.
\end{enumerate}
\end{theorem}
\noindent {\it Proof:} 
Since $\K$ is algebraically closed and since $V$ has finite positive 
dimension, there exists a scalar $\eta \in \K$ such that
$V_{02}(\eta)\not=0$. Observe $\eta \not=0$ since $x_{02}$ is
invertible. By Lemma
\ref{lem:diminv}
we find $V_{01}(\eta)\not=0$. Since $\eta \not=0$ and
since $q$ is not a root of $1$, the scalars $
\eta, q^2\eta, q^4\eta,\ldots $ are mutually distinct,
so they cannot all be eigenvalues for $x_{01}$ on $V$.
Therefore there exists a nonzero $\theta \in \K$ such that
$V_{01}(\theta)\not=0$ and
$V_{01}(q^2\theta)=0$. Similarly there exists an integer
$d\geq 0$ such that $V_{01}(q^{-2n}\theta)$ is nonzero
for $0 \leq n \leq d$ and zero for $n=d+1$. We show that
\begin{eqnarray}
\label{eq:vchain}
V_{01}(\theta)+
V_{01}(q^{-2}\theta)+\cdots +
V_{01}(q^{-2d}\theta)
\end{eqnarray}
is equal to $V$. By
Lemma \ref{lem:qweyl} and since
$V_{01}(q^{-2d-2}\theta)=0$
the space (\ref{eq:vchain})
is invariant under each of $x_{12}, x_{13}$.
The space (\ref{eq:vchain})
is invariant under $x_{31}$ since $x_{31}$ is the
inverse of $x_{13}$.
By Lemma \ref{lem:qweyl2} and since
$V_{01}(q^{2}\theta)=0$
the space (\ref{eq:vchain})
is invariant under each of $x_{20}, x_{30}$.
The space (\ref{eq:vchain})
is invariant under $x_{02}$  since
$x_{02}$ is the inverse of
$x_{20}$.
By Lemma
\ref{lem:key},
 and since each of
$V_{01}(q^2\theta)$,
$V_{01}(q^{-2d-2}\theta)$ is zero,
the space
(\ref{eq:vchain})
is invariant under 
$x_{23}$.
We have now shown that
(\ref{eq:vchain})
is invariant under 
each generator $x_{ij}$ of $\boxtimes_q$, so
(\ref{eq:vchain}) is a $\boxtimes_q$-submodule of $V$.
Each term in
(\ref{eq:vchain}) is nonzero and there is at least one
term so 
(\ref{eq:vchain}) is nonzero. By these comments and
since the $\boxtimes_q$-module $V$ is irreducible
we find 
(\ref{eq:vchain}) is equal to $V$.
This shows that the action of $x_{01}$ on $V$ is semisimple
with eigenvalues
$\Delta:=\lbrace q^{-2n}\theta \;|\;0 \leq n\leq d\rbrace$.
By Lemma
\ref{lem:diminv}, $\Delta$ contains the multiplicative
inverse of each of its elements; therefore $q^{-2d}\theta=\theta^{-1}$
so $\theta^2=q^{2d}$. Consequently there exists a scalar
$\varepsilon \in \lbrace 1,-1\rbrace $ such that
$\theta=\varepsilon q^d$, and we get $\Delta=\lbrace \varepsilon q^{d-2n}\;|\;0 \leq n \leq d\rbrace$. So far we have shown that the action of
$x_{01}$ on $V$ is semisimple with eigenvalues
$\lbrace \varepsilon q^{d-2n}\;|\;0 \leq n \leq d\rbrace$.
Applying Lemma \ref{lem:diminv} we find that for each generator
$x_{ij}$ the action on $V$ is semisimple with
eigenvalues
$\lbrace \varepsilon q^{d-2n}\;|\;0 \leq n \leq d\rbrace$.
\hfill $\Box $  \\

\begin{definition}
\rm
Referring to Theorem
\ref{thm:diminv}, we call $d$ the {\it diameter}
of $V$. We call $\varepsilon$ the {\it type} of
$V$.
\end{definition}

\begin{note}
\label{note1}
\rm
Let $V$denote a finite dimensional irreducible
$\boxtimes_q$-module of type $\varepsilon$. Applying the automorphism
from Lemma
\ref{lem:flip} the type becomes $-\varepsilon$.
\end{note}

\begin{note}
\label{note2}
\rm
In view of Note \ref{note1}, as we proceed we will focus on
the finite dimensional
irreducible $\boxtimes_q$-modules of type 1. 
\end{note}

\section{The shape of a finite dimensional irreducible $\boxtimes_q$-module}

\noindent Let $V$ denote a finite dimensional
 irreducible $\boxtimes_q$-module of type 1.
In this section we associate with $V$ a finite sequence of positive
integers called the {\it shape} of $V$.

\medskip
\noindent We will use the following notation.
Let $V$ denote a vector space over $\K$ with finite positive
dimension. 
Let $(s_0,s_1,\ldots,s_d)$ denote a finite sequence consisting
of positive integers whose sum is the dimension of $V$. 
By a {\it decomposition of $V$ of shape
$(s_0,s_1,\ldots,s_d)$}
we mean a sequence $V_0, V_1,\ldots, V_d$ 
of subspaces of $V$ such that $V_n$ has dimension $s_n$ for
$0 \leq n \leq d$ and
\begin{eqnarray*}
V = V_0+V_1+\cdots +V_d \qquad \qquad (\mbox{\rm direct sum}).
\end{eqnarray*}
We call $d$ the {\it diameter} of the decomposition.
 For
$0 \leq n \leq d$ we call $V_n$ the {\it $n$th component}
 of the decomposition.
For notational convenience we define $V_{-1}=0$, $V_{d+1}=0$.
By the {\it inversion} of $V_0,V_1,\ldots,V_d$
we mean the decomposition
$V_d,V_{d-1},\ldots, V_0$. 

\begin{definition}
\label{def:decij}
\rm
Let $V$ denote a
finite dimensional irreducible $\boxtimes_q$-module of type
$1$ and diameter $d$. For each generator $x_{ij}$ of $\boxtimes_q$
we define a decomposition of $V$ which we call $\lbrack i,j\rbrack$.
The decomposition $\lbrack i,j\rbrack $ has diameter $d$.
For $0 \leq n\leq d$ the $n$th component of
$\lbrack i,j\rbrack$ is the eigenspace of $x_{ij}$ on
$V$ associated with the eigenvalue $q^{d-2n}$. 
\end{definition}

\begin{note}\rm
With reference to Definition
\ref{def:decij}, for $i \in \Z_4$ the decomposition
$\lbrack i,i+2\rbrack $ is the inversion of
$\lbrack i+2,i\rbrack$.
\end{note}

\begin{proposition}
\label{lem:shape}
Let $V$ denote a finite dimensional irreducible $\boxtimes_q$-module
of type 1 and diameter $d$. Choose a generator $x_{ij}$ of $\boxtimes_q$
and consider the corresponding decomposition $\lbrack i,j\rbrack$ of $V$  
from Definition
\ref{def:decij}.
Then the shape of this decomposition is independent
of the choice of generator. Denoting the shape by
$(\rho_0,\rho_1,\ldots,\rho_d)$ we have $\rho_n=\rho_{d-n}$ for $0 \leq n\leq d$.
\end{proposition}
\noindent {\it Proof:}
Immediate from Lemma
\ref{lem:diminv}.
\hfill $\Box $ 

\begin{definition} \rm
\label{def:shapev}
Let $V$ denote a finite dimensional irreducible $\boxtimes_q$-module
of type 1 and diameter $d$. By the {\it shape} of $V$
we mean the sequence $(\rho_0,\rho_1,\ldots,\rho_d)$ from
Proposition
\ref{lem:shape}.
\end{definition}

\section{Finite dimensional irreducible $\boxtimes_q$-modules; the
$x_{ij}$ action}

\noindent Let $V$ denote a finite dimensional irreducible
$\boxtimes_q$-module of type 1. In this section we describe how
 each generator $x_{ij}$ of $\boxtimes_q$ 
 acts on the eigenspaces of the other
generators. We will treat separately the cases $j-i=1$ and $j-i=2$.

\begin{theorem}
\label{thm:sixdecp}
Let $V$ denote a finite dimensional irreducible
$\boxtimes_q$-module of type 1 and diameter $d$.
Let 
$V_0, V_1,\ldots,V_d$ denote a decomposition of
$V$ from Definition
\ref{def:decij}. Then for $i \in \Z_4$ and
for $0 \leq n\leq d$ the action of $x_{i,i+1}$ on
$V_n$ is given as follows.
\medskip

\centerline{
\begin{tabular}[t]{c|c}
       {\rm decomposition} & {\rm action of $x_{i,i+1}$ on $V_n$}
 \\ \hline  \hline
	$\lbrack i,i+1\rbrack$ & $(x_{i,i+1}-q^{d-2n}I)V_n=0$    
	\\
	$\lbrack i+1,i+2\rbrack$ & 
              $(x_{i,i+1}-q^{2n-d}I)V_n \subseteq V_{n-1}$   \\
	$\lbrack i+2,i+3\rbrack$ &
	$x_{i,i+1}V_n\subseteq V_{n-1}+V_n+V_{n+1}$  \\ 
	$\lbrack i+3,i\rbrack$ & 
	$(x_{i,i+1}-q^{2n-d}I)V_n\subseteq V_{n+1}$
	\\ 
        $\lbrack i,i+2\rbrack$ &
	$
	(x_{i,i+1}-q^{d-2n}I)V_n\subseteq V_{n-1}$
	\\ 
	$ \lbrack i+1,i+3\rbrack $ & 
	$
	(x_{i,i+1}-q^{2n-d}I)V_n\subseteq V_{n-1}$
	\end{tabular}}
\medskip
\noindent 
\end{theorem}
\noindent {\it Proof:} We consider each of the six rows
of the table.
\\
\noindent $\lbrack i,i+1\rbrack$: By Definition
\ref{def:decij}
 $V_n$
is the eigenspace for $x_{i,i+1}$ associated with the
eigenvalue $q^{d-2n}$.
\\
\noindent $\lbrack i+1,i+2\rbrack$:
Apply Lemma
\ref{lem:qweyl2} (with $A=x_{i,i+1}$ and $B=x_{i+1,i+2}$).
\\
\noindent $\lbrack i+2,i+3\rbrack$: 
Apply Lemma
\ref{lem:key} (with $A=x_{i+2,i+3}$ and $B=x_{i,i+1}$).
\\
\noindent $\lbrack i+3,i\rbrack$: 
Apply Lemma
\ref{lem:qweyl} (with $A=x_{i+3,i}$ and $B=x_{i,i+1}$).
\\
\noindent $\lbrack i,i+2\rbrack$: 
Apply Lemma
\ref{lem:qweyl} (with $A=x_{i+2,i}$ and $B=x_{i,i+1}$).
\\
\noindent $\lbrack i+1,i+3\rbrack$: 
Apply Lemma
\ref{lem:qweyl2} (with $A=x_{i,i+1}$ and $B=x_{i+1,i+3}$).
\hfill $\Box $ \\

\begin{theorem}
\label{thm:sixdecp2}
Let $V$ denote a finite dimensional irreducible
$\boxtimes_q$-module of type 1 and diameter $d$.
Let 
$V_0, V_1,\ldots,V_d$ denote a decomposition of
$V$ from Definition
\ref{def:decij}. Then for $i \in \Z_4$ and
for $0 \leq n\leq d$ the action of $x_{i,i+2}$ on
$V_n$ is given as follows.
\medskip

\centerline{
\begin{tabular}[t]{c|c}
       {\rm decomposition} & {\rm action of $x_{i,i+2}$ on $V_n$}
 \\ \hline  \hline
	$\lbrack i,i+1\rbrack$ &
$(x_{i,i+2}-q^{d-2n}I)V_n \subseteq V_0+\cdots+V_{n-1}$ 
	\\
	$\lbrack i+1,i+2\rbrack$ & 
       $(x_{i,i+2}-q^{d-2n}I)V_n \subseteq V_{n+1}+\cdots+V_{d}$     \\
	$\lbrack i+2,i+3\rbrack$ &
	$(x_{i,i+2}-q^{2n-d}I)V_n\subseteq V_{n-1}$
         \\ 
	$\lbrack i+3,i\rbrack$ & 
	$(x_{i,i+2}-q^{2n-d}I)V_n\subseteq V_{n+1}$
	\\ 
        $\lbrack i,i+2\rbrack$ &
	$
	(x_{i,i+2}-q^{d-2n}I)V_n=0$
	\\ 
	$ \lbrack i+1,i+3\rbrack $ & 
	$
	x_{i,i+2}V_n\subseteq V_{n-1}+\cdots+V_d$
	\end{tabular}}
\medskip
\noindent 
\end{theorem}
\noindent {\it Proof:} We consider each of the six rows
of the table.
\\
\noindent $\lbrack i,i+1\rbrack$:
By Lemma
\ref{lem:qweyl2} (with $A=x_{i+2,i}$ and $B=x_{i,i+1}$) we
find
\begin{eqnarray*}
(x_{i+2,i}-q^{2r-d}I)V_r\subseteq V_{r-1}
\qquad \qquad 
(0 \leq r\leq d).
\end{eqnarray*}
The result follows from this and since
$x_{i,i+2}$ is the inverse of 
$x_{i+2,i}$.
\\
\noindent $\lbrack i+1,i+2\rbrack$:
By Lemma
\ref{lem:qweyl} (with $A=x_{i+1,i+2}$ and $B=x_{i+2,i}$) we
find
\begin{eqnarray*}
(x_{i+2,i}-q^{2r-d}I)V_r\subseteq V_{r+1}\qquad \qquad 
(0 \leq r\leq d).
\end{eqnarray*}
The result follows from this and since
$x_{i,i+2}$ is the inverse of 
$x_{i+2,i}$.
\\
\noindent $\lbrack i+2,i+3\rbrack$: 
Apply Lemma
\ref{lem:qweyl2} (with $A=x_{i,i+2}$ and $B=x_{i+2,i+3}$).
\\
\noindent $\lbrack i+3,i\rbrack$: 
Apply Lemma
\ref{lem:qweyl} (with $A=x_{i+3,i}$ and $B=x_{i,i+2}$).
\\
\noindent $\lbrack i,i+2\rbrack$: 
Recall that $V_n$ is the eigenspace of $x_{i,i+2}$ associated
with the eigenvalue $q^{d-2n}$.
\\
\noindent $\lbrack i+1,i+3\rbrack$: 
Let $U_0,U_1,\ldots,U_d$ denote the decomposition
$\lbrack i+2,i+3\rbrack$. From row
$\lbrack i+2,i+3\rbrack$ in the table of this theorem we find
\begin{eqnarray}
x_{i,i+2}U_r\subseteq U_{r-1}+U_r \qquad \qquad (0 \leq r\leq d).
\label{eq:c1}
\end{eqnarray}
By Lemma \ref{lem:qweyl3} (with $A=x_{i+2,i+3}$ and $B=x_{i+3,i+1}$),
\begin{eqnarray}
V_r+\cdots+V_d = U_r+\cdots+U_d \qquad \qquad (0 \leq r\leq d).
\label{eq:c2}
\end{eqnarray}
We may now argue
\begin{eqnarray*}
x_{i,i+2}V_n &\subseteq & x_{i,i+2}(V_n+\cdots+V_d)
\\
&=& x_{i,i+2}(U_n+\cdots+U_d)       \qquad \qquad (\mbox{\rm by (\ref{eq:c2})})
\\
&\subseteq& U_{n-1}+\cdots+U_d    \qquad   \qquad \qquad (\mbox{\rm by (\ref{eq:c1})})\\
&=&V_{n-1}+\cdots+V_d   \qquad    \qquad \qquad (\mbox{\rm by (\ref{eq:c2})}).
\end{eqnarray*}
\hfill $\Box $ \\

\section{Flags}

\noindent In this section we recall the notion of a {\it flag}.

\medskip
\noindent Let $V$ denote a vector space over $\K$ with finite
positive dimension. Let $(s_0,s_1,\ldots,s_d)$ denote a sequence
of positive integers whose sum is the dimension of $V$.
By a {\it flag on $V$ of shape $(s_0,s_1,\ldots,s_d)$} we mean
a nested sequence $U_0\subseteq U_1\subseteq \cdots \subseteq
U_d$ of subspaces of $V$ such that the dimension of $U_n$
is $s_0+\cdots+s_n$ for $0 \leq n \leq d$. We call
$U_n$ the $n$th {\it component} of the flag. We call
$d$ the {\it diameter} of the flag. We observe 
$U_d=V$.

\medskip
\noindent The following construction yields a flag on $V$.
Let $V_0,V_1,\ldots,V_d$ denote a decomposition of $V$
of shape $(s_0,s_1,\ldots,s_d)$. Define
\begin{eqnarray*}
U_n=V_0+V_1+\cdots+V_n \qquad \qquad (0 \leq n\leq d).
\end{eqnarray*}
Then the sequence $U_0 \subseteq U_1 \subseteq \cdots \subseteq U_d$
is a flag on $V$ of shape $(s_0,s_1,\ldots, s_d)$.
We say this flag is {\it induced} by the decomposition
 $V_0,V_1,\ldots,V_d$.

\medskip
\noindent We now recall what it means for two flags to be {\it opposite}.
Suppose we are given two flags on $V$ with the same diameter: 
 $U_0 \subseteq U_1 \subseteq \cdots \subseteq U_d$
and 
 $U'_0 \subseteq U'_1 \subseteq \cdots \subseteq U'_d$.
We say these flags are {\it opposite} whenever there exists
a decomposition $V_0,V_1,\ldots, V_d$ of $V$ such that
\begin{eqnarray*}
U_n=V_0+V_1+\cdots+V_n, \qquad \qquad 
U'_n=V_d+V_{d-1}+\cdots+V_{d-n}
\end{eqnarray*}
for $0 \leq n\leq d$. In this case
\begin{eqnarray}
\label{eq:zero}
U_i\cap U'_j = 0 \qquad \mbox{if}\quad i+j<d \qquad \qquad (0 \leq i,j\leq d)
\end{eqnarray}
and
\begin{eqnarray}
\label{eq:recover}
V_n=U_n\cap U'_{d-n} \qquad \qquad (0 \leq n\leq d).
\end{eqnarray}
In particular the decomposition $V_0,V_1,\ldots,V_d$ is uniquely
determined by the given flags.

\section{Flags on finite dimensional irreducible $\boxtimes_q$-modules}

\noindent We  return our attention to the finite
dimensional irreducible $\boxtimes_q$-modules.

\begin{theorem}
\label{thm:fourflags}
Let $V$ denote a finite dimensional irreducible $\boxtimes_q$-module
of type 1 and diameter $d$. Then there exists a collection of
flags on $V$, denoted $\lbrack h \rbrack, h \in \Z_4$, that have 
the following property:
for each generator $x_{ij}$ of $\boxtimes_q$ the decomposition
$\lbrack i,j\rbrack$ of $V$ induces $\lbrack i\rbrack$ and
the inversion of $\lbrack i,j\rbrack$ induces $\lbrack j\rbrack$.
\end{theorem}
\noindent {\it Proof:}
For all $h \in \Z_4$ let $\lbrack h \rbrack$ denote
the flag on $V$ induced by the inversion of
$\lbrack h-1,h\rbrack$. By Lemma
\ref{lem:qweyl3}  (with $A=x_{h-1,h}$ and $B=x_{h,h+1}$) the
flag on $V$ induced by $\lbrack h,h+1\rbrack$ is equal to $\lbrack h \rbrack$.
By Lemma \ref{lem:qweyl3} (with $A=x_{h-1,h}$ and $B=x_{h,h+2}$) the
flag on $V$  induced by
$\lbrack h,h+2\rbrack$ is equal to $\lbrack h \rbrack$.
The result follows.
\hfill $\Box $ \\

\begin{lemma} Let $V$ denote a finite dimensional irreducible
$\boxtimes_q$-module of type 1. Then for $i \in \Z_4$ the shape
of the flag $\lbrack i \rbrack $ coincides with the shape of $V$.
\end{lemma}
\noindent {\it Proof:} 
Let $(\rho_0, \rho_1, \ldots, \rho_d)$ denote
the shape of $V$. By Proposition
\ref{lem:shape} the
decomposition $\lbrack i,i+1 \rbrack$ has shape
$(\rho_0, \rho_1, \ldots, \rho_d)$.
This decompostion
induces the flag $\lbrack i \rbrack$ so
the flag $\lbrack i \rbrack$ 
has 
shape 
$(\rho_0, \rho_1, \ldots, \rho_d)$.
\hfill $\Box $ \\

\begin{theorem} Let $V$ denote a finite
dimensional irreducible $\boxtimes_q$-module of
type 1. Then the flags $\lbrack i \rbrack$, $i \in \Z_4$ on $V$  
from Theorem \ref{thm:fourflags}
are mutually opposite.
\end{theorem}
\noindent {\it Proof:} We invoke Theorem
 \ref{thm:fourflags}. For $i \in \Z_4$ the flags
$\lbrack i\rbrack, \lbrack i+1 \rbrack$ are opposite
since the decomposition $\lbrack i,i+1\rbrack$ induces
$\lbrack i \rbrack$ and the inversion of this decomposition
induces $\lbrack i+1\rbrack$. The flags 
$\lbrack i \rbrack, \lbrack i+2\rbrack$ are opposite
since the decomposition
$\lbrack i,i+2\rbrack$ induces $\lbrack i \rbrack$ and
the inversion of this decomposition induces $\lbrack i+2\rbrack $.
The result follows.
\hfill $\Box $ \\

\begin{theorem}
\label{thm:flagdec}
Let $V$ denote a finite dimensional irreducible $\boxtimes_q$-module
of type 1 and diameter $d$. Pick a generator $x_{ij}$ of $\boxtimes_q$
and consider the corresponding decomposition $\lbrack i,j\rbrack $
of $V$ from Definition \ref{def:decij}.
For $0 \leq n\leq d$ the $n$th component of $\lbrack i,j\rbrack$
is the intersection of the following two sets:
\begin{enumerate}
\item[{\rm (i)}]  
component $n$ of the flag $\lbrack i\rbrack$;
\item[{\rm (ii)}] 
component $d-n$ of the flag $\lbrack j\rbrack$.
\end{enumerate}
\end{theorem}
\noindent {\it Proof:} Combine
Theorem
 \ref{thm:fourflags}
and line 
(\ref{eq:recover}).
\hfill $\Box $ \\

\section{From $\boxtimes_q$-modules to ${\mathcal A}_q$-modules}

In this section we give the proof of
Theorem
\ref{thm:2}. Our proof is based on the following proposition.

\begin{proposition}
\label{prop:w}
Let $V$ denote a finite dimensional irreducible $\boxtimes_q$-module
of type 1. Let $W$ denote a nonzero subspace of $V$ such
that $x_{01}W\subseteq W$ and
$x_{23}W\subseteq W$. Then $W=V$.
\end{proposition}
\noindent {\it Proof:}
Without loss we may assume $W$ is irreducible as a module for $x_{01},x_{23}$.
Let $V_0,V_1,\ldots,V_d$ denote the decomposition
$\lbrack 0,1\rbrack$ and let $V'_0,V'_1,\ldots,V'_d$ denote the
decomposition $\lbrack 2,3\rbrack$. Recall that $x_{01}$ (resp. $x_{23}$)
is semisimple on $V$ with eigenspaces $V_0,V_1,\ldots,V_d$ (resp.
$V'_0,V'_1,\ldots,V'_d$). By this and since $W$ is invariant under
each of $x_{01},x_{23}$ we find
\begin{eqnarray}
\label{eq:wsum}
W=\sum_{n=0}^d W\cap V_n,
\qquad \qquad 
W=\sum_{n=0}^d W\cap V'_n.
\end{eqnarray}
Define
\begin{eqnarray}
\label{eq:m}
m= \mbox{\rm min}\lbrace n \;|\;0 \leq n\leq d, \quad W\cap V_n\not=0\rbrace.
\end{eqnarray}
We claim
\begin{eqnarray}
\label{eq:m1}
m&=& 
\mbox{\rm min}\lbrace n \;|\;0 \leq n\leq d, \quad W\cap V'_n\not=0\rbrace,
\\
\label{eq:m2}
m&=&
 \mbox{\rm min}\lbrace n \;|\;0 \leq n\leq d, \quad W\cap V_{d-n}\not=0
\rbrace,
\\
\label{eq:m3}
m&=&
 \mbox{\rm min}\lbrace n \;|\;0 \leq n\leq d, \quad W\cap V'_{d-n}\not=0
\rbrace.
\end{eqnarray}
To prove the claim we let $m',m'',m'''$ denote the integers on
the right-hand side of 
(\ref{eq:m1})--(\ref{eq:m3}) respectively.
We will show $m,m',m'',m'''$ coincide by
showing 
$m\leq m'\leq m''\leq m'''\leq m$.
Suppose $m>m'$. By
(\ref{eq:m}) and the equation on the left in 
(\ref{eq:wsum}),
the space $W$ is contained
 in component $d-m$ of the flag $\lbrack 1\rbrack$.
By construction $W$ has nonzero intersection with component $m'$
of the flag $\lbrack 2\rbrack$. Since $m>m'$ the component $d-m$ of
$\lbrack 1\rbrack$ has zero intersection with  component
$m'$ of $\lbrack 2\rbrack$, for a contradiction. Therefore
$m\leq m'$.
Next suppose
$m'>m''$. By the definition of $m'$
 and the equation on the right in 
(\ref{eq:wsum}),
the space $W$ is contained
 in component $d-m'$ of the flag $\lbrack 3\rbrack$.
By construction 
 $W$ has nonzero intersection with component $m''$
of the flag $\lbrack 1\rbrack$. Since $m'>m''$ the component $d-m'$ of
$\lbrack 3\rbrack$ has zero intersection with  component
$m''$ of $\lbrack 1\rbrack$, for a contradiction. Therefore
$m'\leq m''$.
Next suppose
$m''>m'''$. By
the definition of $m''$ and the equation on the left in 
(\ref{eq:wsum}),
the space $W$ is contained
 in component $d-m''$ of the flag $\lbrack 0\rbrack$.
By construction $W$ has nonzero intersection with component $m'''$
of the flag $\lbrack 3\rbrack$. Since $m''>m'''$ the component $d-m''$ of
$\lbrack 0\rbrack$ has zero intersection with  component
$m'''$ of $\lbrack 3\rbrack$, for a contradiction. Therefore
$m''\leq m'''$.
Now suppose
$m'''>m$. By
the definition of $m'''$ and the equation on the right in 
(\ref{eq:wsum}),
the space $W$ is contained
 in component $d-m'''$ of the flag $\lbrack 2\rbrack$.
By (\ref{eq:m}) $W$ has nonzero intersection with component $m$
of the flag $\lbrack 0\rbrack$. Since $m'''>m$ the component $d-m'''$ of
$\lbrack 2\rbrack$ has zero intersection with  component
$m$ of $\lbrack 0\rbrack$, for a contradiction. Therefore
$m'''\leq m$.
We have now shown 
$m\leq m'\leq m''\leq m'''\leq m$. Therefore
$m,m',m'',m'''$ coincide 
 and the claim is proved.
The claim implies that for all $i \in \Z_4$ the component
$d-m$ of the flag $\lbrack i \rbrack $ contains $W$,
and component $m$ of $\lbrack i \rbrack$ has nonzero intersection
with $W$.
We can now easily show $W=V$.
Since $V$ is irreducible as a $\boxtimes_q$-module it suffices
to show that $W$ is invariant under $\boxtimes_q$. By construction
$W$ is invariant under each of $x_{01},x_{23}$.
We now let $x_{rs}$ denote one of $x_{12},x_{13},x_{20},x_{30}$
and show $x_{rs}W\subseteq W$. 
Let $W'$ denote the span of the set of vectors in 
$W$ that are eigenvectors for $x_{rs}$. By construction
$W'\subseteq W$ and $x_{rs}W'\subseteq W'$. We show
$W'=W$. To this end we show that
$W'$ is nonzero and invariant under each of $x_{01}, x_{23}$.
We now show $W'\not=0$. By the comment after the preliminary 
claim, $W$ has nonzero intersection with component $m$
of the flag $\lbrack r \rbrack $ and $W$ is contained in
component $d-m$ of the flag $\lbrack s \rbrack$. By 
Theorem
\ref{thm:flagdec} the intersection of component
$m$ of $\lbrack r \rbrack $ and component $d-m$ of
$\lbrack s\rbrack$ is equal to component $m$ of the
decomposition $\lbrack r,s\rbrack$, which is an eigenspace
for $x_{rs}$. The intersection of $W$ with this eigenspace
is nonzero and contained in $W'$, so $W'\not=0$.
We now show $x_{01}W'\subseteq W'$.
To this end we pick $v \in W'$ and show $x_{01}v \in W'$.
Without loss we may assume that $v$ is an eigenvector for $x_{rs}$;
let $\theta$ denote the corresponding eigenvalue. Then
$\theta\not=0$ by Theorem
\ref{thm:diminv}. Recall $v\in W'$ and
 $W'\subseteq W$
so $v \in W$. The space $W$ is $x_{01}$-invariant so $x_{01}v \in W$.
By these comments $(x_{01}-\theta^{-1}I)v \in W$. By Lemma
\ref{lem:qweyl} or Lemma \ref{lem:qweyl2} the vector
$(x_{01}-\theta^{-1}I)v$ is contained in an eigenspace of $x_{rs}$,
so $(x_{01}-\theta^{-1}I)v\in W'$. By this and since
$v \in W'$ we have $x_{01}v \in W'$.
We have now shown $x_{01}W'\subseteq W'$ as desired.
Next we show $x_{23}W'\subseteq W'$.
To this end we pick $u \in W'$ and show $x_{23}u \in W'$.
Without loss we may assume that $u$ is an eigenvector for $x_{rs}$;
let $\eta$ denote the corresponding eigenvalue. Then
$\eta\not=0$ by Theorem
\ref{thm:diminv}. Recall 
$u\in W'$
and $W'\subseteq W$ 
so $u \in W$. The space $W$ is $x_{23}$-invariant so $x_{23}u \in W$.
By these comments $(x_{23}-\eta^{-1}I)u \in W$. By Lemma
\ref{lem:qweyl} or Lemma \ref{lem:qweyl2} the vector
$(x_{23}-\eta^{-1}I)u$ is contained in an eigenspace of $x_{rs}$,
so $(x_{23}-\eta^{-1}I)u\in W'$. By this and since
$u \in W'$ we have $x_{23}u \in W'$.
We have now shown $x_{23}W'\subseteq W'$ as desired.
So far we have shown that
$W'$ is nonzero and invariant under each of
$x_{01},x_{23}$. Now $W'=W$ by the irreducibility of $W$,
so $x_{rs}W\subseteq W$. It follows that
$W$ is $\boxtimes_q$-invariant. The space $V$ is
irreducible as a $\boxtimes_q$-module so $W=V$.
\hfill $\Box $ \\

\noindent It is now a simple matter to prove Theorem
\ref{thm:2}.
\medskip
\noindent

\medskip
\noindent {\it Proof of Theorem 
\ref{thm:2}:} The generators $x_{01},x_{23}$ satisfy the
cubic $q$-Serre relations 
(\ref{eq:qserre}); by this and
(\ref{eq:qs1}), (\ref{eq:qs2}) 
 there exists
an ${\mathcal A}_q$-module structure on
$V$ such that the standard generators $x,y$ act
as $x_{01}$ and $x_{23}$ respectively. This module
structure is unique since $x,y$ generate ${\mathcal A}_q$.
This module structure is irreducible by Proposition
\ref{prop:w}. By Theorem
\ref{thm:diminv}
and the construction, for each of $x,y$ the action on $V$ is semisimple
with eigenvalues $\lbrace q^{d-2n}\;|\;0 \leq n \leq d\rbrace$.
Therefore the ${\mathcal A}_q$-module structure is NonNil
and type $(1,1)$.
\hfill $\Box $ \\

\section{From ${\mathcal A}_q$-modules to $\boxtimes_q$-modules}

\noindent In this section we give the proof of 
Theorem \ref{thm:1}.

\medskip 
\noindent {\it Proof of Theorem 
\ref{thm:1}:}
By \cite[Corollary 2.8]{NN} and since
the ${\mathcal A}_q$-module $V$
is NonNil, the standard generators $x,y$ are
semisimple on $V$.
 Since $V$ has
type $(1,1)$ there exists an integer $d\geq 0$ 
such that for each of $x,y$
 the set of distinct eigenvalues on $V$ is $\lbrace q^{d-2n}\;|\;0 \leq n \leq d\rbrace$.
We now define two decompositions on $V$, denoted
$\lbrack 0,1\rbrack$ and $\lbrack 2,3\rbrack$, 
where we view
$0,1,2,3$ as the elements of $\Z_4$. Each of these decompositions
has diameter $d$. For $0 \leq n\leq d$ the $n$th component of
$\lbrack 0,1\rbrack$ (resp. $\lbrack 2,3\rbrack$) is
the eigenspace for $x$ (resp. $y$) on $V$ 
associated with the eigenvalue $q^{d-2n}$.
We now define some flags on
$V$, denoted $\lbrack i\rbrack$, $i \in \Z_4$.
The flag $\lbrack 0 \rbrack$ is induced by $\lbrack 0,1\rbrack$
and the flag $\lbrack 1\rbrack $ is induced by the inversion of
$\lbrack 0,1\rbrack$. 
The flag $\lbrack 2 \rbrack$ is induced by $\lbrack 2,3\rbrack$
and the flag $\lbrack 3\rbrack $ is induced by the inversion of
$\lbrack 2,3\rbrack$. By construction the flags $\lbrack 0\rbrack,
\lbrack 1\rbrack$ are opposite and the flags
$\lbrack 2\rbrack, \lbrack 3\rbrack$ are opposite.
By \cite[Lemma 4.3]{tdanduq} the flags
$\lbrack i\rbrack$, $i\in \Z_4$ are mutually opposite.
Now for $i,j\in \Z_4$ such that $j-i=1$ or $j-i=2$,
there exists a decomposition  $\lbrack i,j\rbrack$ of $V$
such that $\lbrack i,j\rbrack $ induces $\lbrack i\rbrack$
and the inversion of $\lbrack i,j\rbrack$ induces $\lbrack j\rbrack$.
We note that $\lbrack i,i+2\rbrack$ is the inversion of 
$\lbrack i+2,i\rbrack$ for $i \in \Z_4$.
For $i,j\in \Z_4$ such that $j-i=1$ or $j-i=2$, let $x_{ij}:V\to V$
denote the linear transformation such that for $0 \leq n\leq d$
the $n$th component of $\lbrack i,j\rbrack$ is the eigenspace for
$x_{ij}$ associated with the eigenvalue $q^{d-2n}$. We observe
that the standard generators $x,y$ act on $V$ as $x_{01}$
and $x_{23}$ respectively. We now show that the above 
transformations $x_{ij}$
satisfy the defining relations
for $\boxtimes_q$. The  $x_{ij}$
satisfy Definition
\ref{def:qtet}(i) by the construction.
The  $x_{ij}$ satisfy
Definition
\ref{def:qtet}(ii) by \cite[Theorems 7.1,10.1,10.2]{tdanduq}.
The $x_{ij}$ satisfy
Definition
\ref{def:qtet}(iii) by
\cite[Theorem 12.1]{tdanduq}.
We have now shown that the transformations
$x_{ij}$ satisfy the defining relations
for $\boxtimes_q$. Therefore they induce
 a $\boxtimes_q$-module
structure on $V$.
So far we have shown that there exists a $\boxtimes_q$-module
structure on $V$ such that the standard generators 
$x$ and $y$
act as $x_{01}$ and $x_{23}$ respectively.
Next we show that this $\boxtimes_q$-module structure
is unique. Suppose we are given
any $\boxtimes_q$-module structure on $
V$ such that the standard generators
$x$ and $y$
act as $x_{01}$ and $x_{23}$ respectively.
This $\boxtimes_q$-module structure is irreducible by construction
and since the 
${\mathcal A}_q$-module $V$ is irreducible.
This $\boxtimes_q$-module structure is 
type $1$ and diameter $d$, since the action
of $x_{01}$ on $V$ has eigenvalues $\lbrace q^{d-2n} \;|\;0 \leq n \leq d
\rbrace $. For each generator $x_{ij}$ of $\boxtimes_q$ the
action on $V$ is determined by the decomposition $\lbrack i,j\rbrack$.
By Theorem
\ref{thm:flagdec} the decomposition $\lbrack i,j\rbrack$
is determined by the flags  $\lbrack i \rbrack $ and
$\lbrack j \rbrack$.
Therefore our $\boxtimes_q$-module structure on $V$
is determined by the flags $\lbrack h\rbrack, h \in \Z_4$.
By construction the flags $\lbrack 0\rbrack$ and $\lbrack 1 \rbrack$
are determined by the decomposition $\lbrack 0,1 \rbrack $
and hence by the action of $x$ on $V$.
Similarly 
 the flags $\lbrack 2\rbrack$ and $\lbrack 3 \rbrack$
are determined by the decomposition $\lbrack 2,3 \rbrack $
and hence by the action of $y$ on $V$.
Therefore the given $\boxtimes_q$-module structure on $V$
is determined by the action of $x$ and $y$ on $V$,
so this 
$\boxtimes_q$-module structure is unique.
We have now shown that 
there exists
a unique
 $\boxtimes_q$-module
structure on $V$ such that the standard generators 
$x$ and $y$
act as $x_{01}$ and $x_{23}$ respectively.
We mentioned earlier that this
 $\boxtimes_q$-module
structure is irreducible and has type $1$.
\hfill $\Box $ \\

\section{Suggestions for further research}

\noindent In this section we give some suggestions for further
research.

\begin{problem}
\rm
Find a basis for the $\K$-vector space $\boxtimes_q$.
\end{problem}

\begin{problem} \rm
Find all the Hopf-algebra structures on $\boxtimes_q$.
\end{problem}

\begin{problem} \rm Find the automorphism group of
the $\K$-algebra $\boxtimes_q$.
\end{problem}

\begin{problem} \rm Find all the 2-sided ideals 
of the $\K$-algebra 
$\boxtimes_q$.
\end{problem}

 \begin{conjecture}
\rm Let $V$ denote a finite dimensional irreducible $\boxtimes_q$-module
of type $1$. Then there exists a nondegenerate symmetric bilinear
form $(\,,\,)$ on $V$ such that
\begin{eqnarray*}
(r.u,v)=(u,\omega(r).v) \qquad \qquad r\in \boxtimes_q, \quad u,v \in V,
\end{eqnarray*}
where $\omega$ is the antiautomorphism from Lemma
\ref{lem:omega}.
\end{conjecture}

\begin{conjecture} \rm 
Let $\Omega$ (resp. $\Omega'$) (resp. $\Omega''$)
denote the subalgebra of $\boxtimes_q$ generated
by $x_{01}$ and $x_{23}$
(resp.
by $x_{12}$ and $x_{30}$)
(resp.
by $x_{02}, x_{20}, x_{13}, x_{31}$).
Then the map
\begin{eqnarray*}
\Omega \otimes
\Omega' \otimes
\Omega'' \quad  &\to & \quad \; \boxtimes_q \\
u\otimes v \otimes w \quad \;\;&\to& \quad uvw 
\end{eqnarray*}
is an isomorphism of $\K$-vector spaces.
By $\otimes$ we mean $\otimes_{\K}$.
\end{conjecture}

\begin{problem} \rm
For $i \in \Z_4$ let $\boxtimes_q^{(i)}$ denote the
set of elements $r \in \boxtimes_q$ that have the
following property: for all finite dimensional irreducible
$\boxtimes_q$-modules $V$ of type 1, $r$ leaves invariant
each component of  $\lbrack i \rbrack $, where $\lbrack i \rbrack$
is the flag on $V$
from
Theorem
\ref{thm:fourflags}.
We observe that 
$\boxtimes_q^{(i)}$
is a subalgebra of $\boxtimes_q$ 
that contains each of the generators
$x_{i,i+1}, x_{i,i+2}, x_{i+2,i}, x_{i-1,i}$.
Find a generating set for the $\K$-algebra
$\boxtimes_q^{(i)}$.
Find a basis for the $\K$-vector space 
$\boxtimes_q^{(i)}$.
For each subset $S \subseteq \Z_4$ 
describe the subalgebra $\cap_{i \in S} \boxtimes_q^{(i)}$.
\end{problem}

\begin{problem} \rm Recall that for each element $r \in \boxtimes_q$
the {\it centralizer} $C(r)$ is the subalgebra of $\boxtimes_q$
consisting of the elements in $\boxtimes_q$ that commute with $r$. For each 
generator $x_{ij}$ of $\boxtimes_q$ describe $C(x_{ij})$. Find
a generating set for the $\K$-algebra
$C(x_{ij})$. Find a 
basis for the $\K$-vector space 
$C(x_{ij})$. 
\end{problem}

\begin{problem} \rm Recall that an element $r \in\boxtimes_q$ is
{\it central} whenever it commutes with each element of $\boxtimes_q$.
The {\it center} of $\boxtimes_q$ is the subalgebra of
$\boxtimes_q$ consisting of the central elements. Find 
the center of $\boxtimes_q$.
\end{problem}

\noindent Before we state our last problem we have two remarks.

\begin{remark} \rm
\label{remark1} 
Let $V$ denote a finite dimensional irreducible $\boxtimes_q$-module.
Then the dual vector space $V^*$ has a $\boxtimes_q$-module
structure such that
\begin{eqnarray*}
(r.f)(v) = f(\omega(r).v) \qquad \qquad r \in \boxtimes_q, 
\quad f \in V^*, \quad v \in V,
\end{eqnarray*}
where $\omega$ is the antiautomorphism of $\boxtimes_q$ from
Lemma
\ref{lem:omega}.
\end{remark}

\begin{remark} \rm
\label{remark2} 
Let $V$ denote a finite dimensional irreducible
$\boxtimes_q$-module.  For each automorphism
$\sigma $ of $\boxtimes_q$ there exists
a $\boxtimes_q$-module structure on $V$, called {\it $V$ twisted via
$\sigma$}, that has the following property: for all $r \in \boxtimes_q$
and for all $v \in V$ the vector $r.v$ computed in $V$ twisted via
$\sigma$ coincides with $\sigma(r).v$ computed in the original
$\boxtimes_q$-module $V$.
\end{remark}

\begin{problem} \rm
Let $V$ denote a finite dimensional irreducible $\boxtimes_q$-module
of type 1, and let the automorphism 
$\rho$ of $\boxtimes_q$ be as in Lemma
\ref{lem:rho}. Using
Remark
\ref{remark1} 
and
Remark
\ref{remark2} 
 we obtain eight $\boxtimes_q$-module 
structures on $V$; these are 
$V$ twisted via $\rho^n$ for $0 \leq n\leq 3$ 
and $V^*$ twisted via $\rho^n$ for $0 \leq n\leq 3$.
How are these
eight $\boxtimes_q$-module structures related up to isomorphism?
\end{problem}

\noindent Tatsuro Ito \hfil\break
\noindent Department of Computational Science \hfil\break
\noindent Faculty of Science \hfil\break
\noindent Kanazawa University \hfil\break
\noindent Kakuma-machi \hfil\break
\noindent Kanazawa 920-1192, Japan \hfil\break
\noindent email:  {\tt ito@kappa.s.kanazawa-u.ac.jp}

\bigskip

\noindent Paul Terwilliger \hfil\break
\noindent Department of Mathematics \hfil\break
\noindent University of Wisconsin \hfil\break
\noindent Van Vleck Hall \hfil\break
\noindent 480 Lincoln Drive \hfil\break
\noindent Madison, WI 53706-1388 USA \hfil\break
\noindent email: {\tt terwilli@math.wisc.edu }\hfil\break

\end{document}